\documentclass[a4paper,12pt]{amsart}

\usepackage{graphicx}

\newtheorem{theorem}{Theorem}[section]

\newtheorem{lemma}[theorem]{Lemma}
\newtheorem{conjecture}[theorem]{Conjecture}
\newtheorem{approximation}[theorem]{Approximation}

\newtheorem{Conclusion}[theorem]{Conclusion}

\newtheorem{Example}[theorem]{Example}
\newenvironment{example}{\begin{Example}\rm}{\end{Example}}
\newtheorem{Remark}[theorem]{Remark}

\newtheorem{Remarks}[theorem]{Remarks}

\newtheorem{Question}[theorem]{Question}

\newcommand{\al}{\alpha}
\newcommand{\be}{\beta}
\newcommand{\ga}{\gamma}
\newcommand{\Ga}{\Gamma}
\newcommand{\de}{\delta}
\newcommand{\De}{\Delta}
\newcommand{\eps}{\varepsilon}

\newcommand{\la}{\lambda}
\newcommand{\La}{\Lambda}
\newcommand{\om}{\omega}
\newcommand{\Om}{\Omega}
\newcommand{\si}{\sigma}
\newcommand{\Si}{\Sigma}
\newcommand{\ze}{\zeta}

\let\cal=\mathcal
\let\Bbb=\mathbb

\begin{document}

\title[Average prime-pair counting formula, 25.02.2009]{
Average prime-pair counting formula}

\date{February 25, 2009}

\author{Jaap Korevaar and Herman te Riele}



\begin{abstract}
Taking $r>0$, let $\pi_{2r}(x)$ denote the number of prime
pairs $(p,\,p+2r)$ with $p\le x$. The prime-pair
conjecture of Hardy and Littlewood (1923) asserts
that $\pi_{2r}(x)\sim 2C_{2r}\,{\rm li}_2(x)$ with an
explicit constant $C_{2r}>0$. There seems to be no good
conjecture for the remainders $\om_{2r}(x)=\pi_{2r}(x)-
2C_{2r}\,{\rm li}_2(x)$ that corresponds to Riemann's formula for
$\pi(x)-{\rm li}(x)$. However, there is a heuristic approximate
formula for averages of the remainders $\om_{2r}(x)$ which is
supported by numerical results.
\end{abstract}

\maketitle

\setcounter{equation}{0}    
\section{Introduction} \label{sec:1}
For $r\in{\Bbb N}$, let $\pi_{2r}(x)$ denote the number of
prime pairs $(p,\,p+2r)$ with $p\le x$. The famous
prime-pair conjecture (PPC) of Hardy and
Littlewood \cite{HL23} asserts that for $x\to\infty$,
\begin{equation} \label{eq:1.1}
\pi_{2r}(x)\sim 2C_{2r}{\rm li}_2(x)=
2C_{2r}\int_2^x\frac{dt}{\log^2 t}\sim
2C_{2r}\frac{x}{\log^2 x}.
\end{equation}
Here $C_2$ is the `twin-prime constant',
\begin{equation} \label{eq:1.2}
C_2 = \prod_{p\,{\rm
prime},\,p>2}\,\left\{1-\frac{1}{(p-1)^2}\right\}
\approx 0.6601618158,
\end{equation}
and the general `prime-pair constant' $C_{2r}$ is given by
\begin{equation} \label{eq:1.3}
C_{2r} = C_2\prod_{p\,{\rm prime},\,p|r,\,p>2}\,\frac{p-1}{p-2}.
\end{equation}

Assuming that the PPC is true, let $\om_{2r}(x)$ denote the remainder
\begin{equation} \label{eq:1.4}
\om_{2r}(x)=\pi_{2r}(x)-2C_{2r}{\rm li}_2(x).
\end{equation}
We have not been able to find a good approximation of the remainders
$\om_{2r}(x)$ that correspond to Riemann's approximate formula
for $\pi(x)-{\rm li}(x)$ (see (\ref{eq:1.9}) below).
Instead, by complex analysis and heuristic arguments we obtain
the following {\it plausible approximation} for {\it averages}
$(1/N)\sum_{r=1}^N\,\om_{2r}(x)$ with large $N$, and we support the
formula by extensive numerical results.
\begin{approximation} \label{approx:1.1}
For $N\ge 1$ and $x\ge N^{2+\de}$, with $0<\de\le 1$, one has
\begin{align} & \label{eq:1.5}
\frac{1}{N}\sum_{r=1}^N\,\{\pi_{2r}(x)-2C_{2r}{\rm li}_2(x)\} =
-\{4+\cal{O}(N^{-1/2}\log x)\}\sum_\rho\,\rho\,{\rm li}_2(x^\rho)
\notag \\ & \qquad\qquad\qquad\qquad 
-\{1+\cal{O}(N^{-1/2}\log x)\}{\rm li}_2(x^{1/2})+\cal{O}(x^{1/(2+\de)}),
\end{align}
with a symmetric sum over the complex zeros $\rho$ of $\ze(s)$.
\end{approximation}
To test this conjectured approximation we observe that
\begin{align} \label{eq:1.6}
\frac{\sum_\rho\,\rho\,{\rm li}_2(x^\rho)}{{\rm li}_2(x^{1/2})}
&= \bigg\{\frac{1}{4}+\cal{O}\bigg(\frac{1}{\log x}\bigg)\bigg\}T(x)
+\cal{O}\bigg(\frac{1}{\log x}\bigg),\;\;\mbox{where} \notag \\
T(x) &= \sum_\rho\,\frac{x^{\rho-1/2}}{\rho}.
\end{align}
Neglecting the $\cal{O}$-terms in (\ref{eq:1.5}) and (\ref{eq:1.6}), 
dividing by ${\rm li}_2(x^{1/2})$, 
and adding $T(x)+1$, we obtain the error function
\begin{equation}\label{eq:1.7}
\De_N(x) \stackrel{\mathrm{def}}{=} \frac{\sum_{r=1}^N\,\om_{2r}(x)}{N{\rm
li}_2(x^{1/2})} +T(x)+1.
\end{equation}
We have evaluated and plotted this function for fixed $x=10^6,\,10^8,\,10^{10},\,10^{12}$ and $2\le 2N\le 5000$
(Figures 1--4 in Section 7), and for fixed $N=400,\,2500$ and $6\le\log_{10}x\le12$
(Figures 5--6 in Section 7).
Taking into account the $\cal{O}$-terms in (\ref{eq:1.5}) and (\ref{eq:1.6}),
$\De_N(x)$ should have the form
$$
\cal{O}\big(N^{-1/2}(\log x)+1/\log x\big)T(x)
+\cal{O}\big(N^{-1/2}(\log x)+1/\log x\big)+\cal{O}(x^{-\de/7}).
$$
Our four plots for fixed $x$, and the two for fixed $N$, show that
Approximation \ref{approx:1.1} is good for large $N$, provided $x/N^2$ is large.

When $x$ is comparable to $N^2$, the theory predicts sizeable
deviation, roughly 
\begin{equation}\label{eq:1.8}
\De_N(x)\approx \overline\De_N(x) \stackrel{\mathrm{def}}{=}
-\frac{2N\log^2 x}{8x^{1/2}\log^2 2N};
\end{equation}
see Section \ref{sec:5}.
Behavior of this type is seen in the plots for $x=10^6$ and $10^8$. 

In connection with (\ref{eq:1.5}) we recall Riemann's
approximation for the remainder $\pi(x)-{\rm li}(x)$. If ${\rm
Re}\,\rho=1/2$ for all
$\rho$ one has
\begin{equation} \label{eq:1.9}
\om(x) = \pi(x)-{\rm li}(x)= -\sum_\rho\,{\rm li}(x^\rho)
-(1/2)\,{\rm li}(x^{1/2})+\cal{O}(x^b),
\end{equation} 
where $b$ may be any number greater than $1/3$. This
follows from von Mangoldt's formula for $\psi(x)=\sum_{n\le
x}\,\La(n)$. Here $\La(n)$ denotes his function: $\La(n)=\log p$
if $n=p^\al$ with $p$ prime, and $\La(n)=0$ if $n$ is not a prime
power. Von Mangoldt proved that
\begin{equation} \label{eq:1.10}
\psi(x)= x -\sum_\rho\,\frac{x^\rho}{\rho}
-\frac{\ze'(0)}{\ze(0)}+\sum_k\frac{x^{-2k}}{2k},
\end{equation}
and this is exact for all $x>1$ where $\psi(x)$ is continuous;
cf.\ Davenport \cite{Da00}, Edwards \cite{Ed74}, Ivi\'{c}
\cite{Iv03}.

\medskip
\begin{center}
{\bf PART I. HEURISTICS}
\end{center}

\setcounter{equation}{0}    
\section{First step towards Conjectured Approximation \ref{approx:1.1}}
\label{sec:2} 
Let us start by introducing the functions
\begin{align} \label{eq:2.1}
\psi_{2r}(x) &= \sum_{n\le x}\,\La(n)\La(n+2r),
\quad\theta_{2r}(x) =
\sum_{p,\,p+2r\,{\rm prime};\;p\le x}\,\log^2 p,\notag \\ 
\theta^*_{2r}(x) &= \sum_{p,\,
p^2\pm 2r\,{\rm prime};\; p\le x}\,\log^2 p.
\end{align}
Partial summation or integration by parts shows that the PPC
(\ref{eq:1.1}) is equivalent to each of the asymptotic relations 
\begin{equation} \label{eq:2.2}
\theta_{2r}(x) \sim 2C_{2r}x,\quad\psi_{2r}(x)\sim
2C_{2r}x\quad\mbox{as}\;\;x\to\infty.
\end{equation}
 
We have counted the prime pairs $(p,\,p+2r)$ with $2r\le 5\cdot 10^3$ and
$p\le x=10^3,\,10^4,\,\cdots,\,10^{12}$. Table 1 is based on
this work; cf.\ also a table in Granville and Martin
\cite{GM06} and one by Fokko van de Bult \cite{Bu07}. The bottom
line shows (rounded) values $L_2(x)$ of the comparison function
$2C_2{\rm li}_2(x)$ for $\pi_2(x)$. 
Computations based on these prime-pair counts make
it plausible that for every $r\in{\Bbb N}$ and every $\eps>0$,
\begin{equation} \label{eq:2.3}
\om_{2r}(x)=\pi_{2r}(x)-2C_{2r}{\rm li}_2(x) = \cal{O}(x^{(1/2)+\eps}).
\end{equation}
Equivalently, there would be inequalities
\begin{equation} \label{eq:2.4}
\Om_{2r}(x)\stackrel{\mathrm{def}}{=}\psi_{2r}(x)-2C_{2r}x =\cal{O}(x^{(1/2)+\eps'}),
\end{equation}
which follow from similar estimates for
$\theta_{2r}(x)-2C_{2r}x$.

Our work requires a good estimate for the difference
$\psi_{2r}(x)-\theta_{2r}(x)$. The main contributions come from the
prime pairs $(p,\,p^2\pm 2r)$ and the sum
$$
\sum_{p,\,p+2r\,{\rm prime};\; p\le x}\,\{\log p\,\log (p+2r)-\log^2 p\}
=\int_2^x \log t\,\log(1+2r/t)\,d\pi_{2r}(t).
$$
Taking $x>2r$ and writing $\int_2^x=\int_{2r}^x+\int_2^{2r}$, one finds
that uniformly in $r$, 
\begin{align} \label{eq:2.5}
\psi_{2r}(x) -\theta_{2r}(x) &= 2\theta^*_{2r}(x^{1/2})
+4rC_{2r}\log\log x \notag \\ &\quad
+\cal{O}(rC_{2r}/\log 2r)+\cal{O}(x^{1/3}\log^2 x).
\end{align}
\begin{table} \label{table:1}
\begin{tabular}{rrrrrrrl} 
$2r\backslash x$ & $10^3$ & $10^4$ & $10^6$ & $10^8$ & $10^{10}$ & $10^{12}$
 & $C_{2r}/C_2$ \\
      &     &       &         &         &            &            &    \\
2     & 35  & 205   &  8169   & 440312  & 27412679   & 1870585220 & 1  \\ 
4     & 41  & 203   &  8144   & 440258  & 27409999   & 1870585459 & 1  \\ 
6     & 74  & 411   & 16386   & 879908  & 54818296   & 3741217498 & 2   \\
8     & 38  & 208   &  8242   & 439908  & 27411508   & 1870580394 & 1   \\
10    & 51  & 270   & 10934   & 586811  & 36548839   & 2494056601 & 4/3 \\
12    & 70  & 404   & 16378   & 880196  & 54822710   & 3741051790 & 2   \\ 
14    & 48  & 245   &  9878   & 528095  & 32891699   & 2244614812 & 6/5 \\
16    & 39  & 200   &  8210   & 441055  & 27414828   & 1870557044 & 1   \\
18    & 74  & 417   & 16451   & 880444  & 54823059   & 3741063106 & 2   \\
20    & 48  & 269   & 10972   & 586267  & 36548155   & 2494072774 & 4/3  \\
22    & 41  & 226   &  9171   & 489085  & 30459489   & 2078443752 & 10/9 \\
24    & 79  & 404   & 16343   & 880927  & 54823858   & 3741122743 & 2    \\ 
30    & 99  & 536   & 21990   &1173934  & 73094856   & 4988150875 & 8/3  \\ 
210   &107  & 641   & 26178   &1409150  & 87712009   & 5985825351 & 16/5 \\  
      &     &       &         &         &            &            & \\
$L_2(x)\,$: & 46 & 214 & 8248  & 440368  & 27411417  & 1870559867 & \\
       &     &       &        &        &        &     &  
\end{tabular}
\caption{$\pi_{2r}(x)$ for selected values of $2r$ and $x$}
\end{table} 

Our first goal will be to motivate the following
\begin{approximation} \label{approx:2.1}
Taking $N$ large and $x$ much larger than $N$, one has
\begin{align} & \label{eq:2.6}
\frac{1}{N}\sum_{r=1}^N\,\{\psi_{2r}(x)-2C_{2r}\,x\} =
-\{4+\cal{O}(N^{-1/2}\log x)\}\sum_\rho\,\frac{x^\rho}{\rho}
\notag \\ &
\quad + \cal{O}(N^{-1/2}x^{1/2}\log x) - \{1+o(1)\}N.
\end{align}
\end{approximation}
 
In support of this conjectured approximation we will derive a related
conjecture involving Dirichlet series. For $s=\si+i\tau$ with
$\si>1/2$, set 
\begin{equation} \label{eq:2.7}
D_{2r}(s) = \sum_{n=1}^\infty\,\frac{\La(n)\La(n+2r)}{n^{2s}}
= \int_1^\infty x^{-2s}d\psi_{2r}(x) 
= 2s\int_1^\infty x^{-2s-1}\psi_{2r}(x)dx.
\end{equation}
Here we use denominator $n^{2s}$ (and not $n^s$) because of the
function $\Phi^\la(s)$ in Theorem \ref{the:3.1} and the
corresponding integral in (\ref{eq:8.2}). 

By a two-way Wiener--Ikehara theorem for Dirichlet series with
positive coefficients, the PPC in the form (\ref{eq:2.2}) is true
if and only if the difference
\begin{equation} \label{eq:2.8}
G_{2r}(s)\stackrel{\mathrm{def}}{=} D_{2r}(s)-\frac{2C_{2r}}{2s-1}
= \int_1^\infty x^{-2s}d\Om_{2r}(x)
= 2s\int_1^\infty x^{-2s-1}\Om_{2r}(x)dx
\end{equation} 
has `good' boundary behavior as $\si\searrow 1/2$. That
is, $G_{2r}(\si+i\tau)$ should tend to a distribution
$G_{2r}\{(1/2)+i\tau\}$ which is locally equal to a
pseudofunction; see Korevaar \cite{Ko05}. Here, a pseudofunction
is the distributional Fourier transform of a bounded function
which tends to zero at infinity. It
cannot have poles and is locally given by Fourier series
whose coefficients tend to zero. In
particular $D_{2r}(s)$ itself would have to show
pole-type behavior, with residue $C_{2r}$, for angular
approach of $s$ to $1/2$ from the right; there should be
no other poles on the line $\{\si=1/2\}$. 

In view of the expected estimate (\ref{eq:2.4}) it
is reasonable to suppose that the difference $G_{2r}(s)$ is
actually analytic for $\si>1/4$. Where would one expect the first
singularities? Assuming Riemann's Hypothesis (RH), we will
motivate a conjecture involving {\it averages} of functions
$G_{2r}(s)$:
\begin{conjecture} \label{con:2.2}
For $\si>1/4$ and  $N\to \infty$, one has
\begin{align} \label{eq:2.9}
\frac{1}{N}\sum_{r=1}^N\,G_{2r}(s) &= \frac{\cal{O}(N^{-1/2})}{(4s-1)^2}
+\frac{\cal{O}(N^{-1/2})}{4s-1} +
\cal{O}(N^{-1/2})\sum_{\rho}\,\frac{1}{(2s-\rho)^2} \notag \\ & \quad
-\{4+\cal{O}(N^{-1/2})\}\,\sum_\rho\,\frac{1}{2s-\rho}+H^N(s),
\end{align}
with symmetric sums over zeta's complex zeros $\rho$. The remainder
$H^N(s)$ has `good' boundary behavior as $\si\searrow 1/4$. For large
$N$ its most significant part may be a term $-\{1+o(1)\}N$.
\end{conjecture}
This conjecture motivates Approximation \ref{approx:2.1}
through formal Fourier inversion; cf.\ (\ref{eq:2.8}). If $L(c)$ denotes
a `vertical line' given by $\si=c>1/4$, then 
\begin{equation} \label{eq:2.10}
\Om_{2r}(x)=\frac{1}{2\pi
i}\int_{L(c)}\,G_{2r}(s)x^{2s}\,\frac{ds}{s}.
\end{equation}

\setcounter{equation}{0} 
\section{The theorem behind Conjecture \ref{con:2.2}} 
\label{sec:3}
To arrive at (\ref{eq:2.9}) we start with a result for a weighted
sum of functions $D_{2r}(s)$; cf.\ \cite{Ko07} (where there is
a less precise result) and the Appendix. The weights are derived
from an even `sieving function'
$E(\nu)$, with $E(0)=1$ and support $[-1,1]$, that can be made to
approach $1$ on $(-1,1)$. A minimal smoothness requirement is
that $E(\nu)$ be absolutely continuous, with derivative $E'(\nu)$
of bounded variation.
\begin{theorem} \label{the:3.1}
Assume RH. Then for $\la>0$ and $1/2<\si<1$,
\begin{align} \label{eq:3.1}
\Phi^\la(s) & \stackrel{\mathrm{def}}{=} D_0(s)+
2\sum_{0<2r\le\la}\,E(2r/\la)D_{2r}(s) \notag \\ &=
\frac{2A^E\la}{2s-1}-4A^E\la \sum_\rho\,
\frac{1}{2s-\rho}+\Si^\la(s)+H^\la_0(s).
\end{align}
The function $D_0(s)$ is obtained from {\rm (\ref{eq:2.7})} by
taking $r=0$ and the constant $A^E$ is given by $\int_0^1
E(\nu)d\nu$. The function $\Si^\la(s)=\Si^{\la,E}(s)$ 
will be described below, and the remainder
$H^\la_0(s)=H^{\la,E}_0(s)$ is analytic for $0<\si<1$. For large $\la$
its most significant part may be a term that behaves like $-\la^2/2$ when
$E(\nu)$ is close to $1$ on $(-1,1)$.
\end{theorem}
The function $D_0(s)$ can be written as follows:
\begin{equation} \label{eq:3.2}
D_0(s)=\sum_{n=1}^\infty\,\frac{\La^2(n)}{n^{2s}}
=\frac{1}{2}\,\frac{d}{ds}
\bigg\{\frac{\ze'(2s)}{\ze(2s)}-\frac{1}{2}\,
\frac{\ze'(4s)}{\ze(4s)}\bigg\}+H_1(s),
\end{equation}
where $H_1(s)$ is analytic for $\si>1/6$. Hence
$D_0(s)$ is meromorphic for $\si>1/6$. Its poles there are purely
quadratic, and located at $s=1/2,\,1/4$ and the points
$\rho/2$. Thus by (\ref{eq:3.1}), and under assumption
(\ref{eq:2.4}), the pole of the difference 
$\Si^\la(s)-D_0(s)$ at $s=1/2$ can only be of first order. The
residue will be given by
\begin{equation} \label{eq:3.3}
R^E(1/2,\la)=2\sum_{0<2r\le\la}\, E(2r/\la)C_{2r}-A^E\la.
\end{equation}
We need the important fact that the constants $C_{2r}$ have
mean value one. Stronger results were obtained by 
Bombieri--Davenport \cite{BoD66} and Montgomery \cite{Mo71}, and these
were later improved by Friedlander and Goldston \cite{FG95} to
\begin{equation} \label{eq:3.4} 
S_m=\sum_{r=1}^m C_{2r}=m-(1/2)\log m+\cal{O}\{\log^{2/3}
(m+1)\}.
\end{equation}
Partial summation in (\ref{eq:3.3}) will thus show that for our
sieving functions $E$, the residue $R(1/2,\la)$ is
$o(\la)$ and in fact, $\cal{O}(\log \la)$ as $\la\to\infty$. We
will see in Section \ref{sec:4} that the outcome $o(\la)$ -- or
$\cal{O}(\la^\eps)$ -- {\it would follow heuristically} from the way
the parameter $\la$ occurs in the terms of $\Si^\la(s)$.

The description of $\Si^\la(s)$ requires a {\it Mellin
transform} associated with the Fourier transform $\hat E^\la(t)$
of $E(\nu/\la)$. For $z=x+iy$ with $0<x<1$ we set
\begin{align} \label{eq:3.5}
M^\la(z) &\stackrel{\mathrm{def}}{=}
\frac{1}{\pi}\int_0^\infty \hat E^\la(t)t^{-z}dt \notag \\ &=
\frac{2}{\pi}\la^z\Ga(1-z)\sin(\pi z/2)\int_0^1E(\nu)\nu^{z-1}dz \\ &=
\frac{2}{\pi}\la^z\Ga(-z-1)\sin(\pi
z/2)\int_0^{1+}\nu^{z+1}dE'(\nu).\notag
\end{align}
The function $M^\la(z)$ extends to a
meromorphic function for $x>-1$ with simple poles at the
points $z=1,\,3,\,\cdots$. The residue of the pole at 
$z=1$ is $-2(\la/\pi)A^E$ with $A^E=\int_0^1 E(\nu)d\nu$,
and $M^\la(0)=1$. Furthermore, the standard order
estimates
\begin{equation} \label{eq:3.6}
\Ga(z)\ll |y|^{x-1/2}e^{-\pi|y|/2},\quad\sin(\pi z/2)\ll
e^{\pi|y|/2}
\end{equation}
for $|x|\le C$ and $|y|\ge 1$ imply the useful
majorization
\begin{equation} \label{eq:3.7}
M^\la(x+iy)\ll \la^x(|y|+1)^{-x-3/2}
\quad\mbox{for}\;\;-1<x\le
C,\;\; |y|\ge 1.
\end{equation}
\begin{example} \label{exam:3.2}
One may take $E(\nu/\la)$ equal to the Fej\'{e}r kernel for ${\Bbb
R}$:
$$ 
E(\nu/\la) = 
\frac{1}{\pi}\int_0^\infty\frac{\sin^2(\la t/2)}
{\la(t/2)^2}\cos \nu t\,dt
=\left\{\begin{array}{ll}
1-|\nu|/\la & \mbox{for $|\nu|\le\la$},
\\ 0 & \mbox{for $|\nu|\ge\la$.}
\end{array}\right.
$$
In this case one finds
$$
M^\la(z) = 
\frac{2}{\pi}\la^z\Ga(-z-1)\sin(\pi z/2).
$$
\end{example}

\noindent {\it The function} $\Si^\la(s)$.
For any $\la>0$, the function $\Si^\la(s)$ is given by the sum
\begin{align} & \label{eq:3.8}
\left\{\frac{\ze'(s)}{\ze(s)}\right\}^2 
+2\,\frac{\ze'(s)}{\ze(s)}\,
\sum_{\rho}\,\Ga(\rho-s)M^\la(\rho-s)
\cos\{\pi(\rho-s)/2\} \notag \\ & 
+ \sum_{\rho',\,\rho''}\,\Ga(\rho'-s)\Ga(\rho''-s) 
M^\la(\rho'+\rho''-2s)\cos\{\pi(\rho'-\rho'')/2\}.
\end{align}
Here $\rho$, $\rho'$ and $\rho''$ independently run over
the complex zeros of $\ze(s)$. It is convenient to denote
the sum of the first two terms by $\Si^\la_1(s)$; for
$0<\si\le 1$ it has poles at $s=1$ and the points
$\rho$. For well-behaved functions $M^\la(z)$, 
the double series defines a function $\Si^\la_2(s)$ as a 
limit of square partial sums. Under RH the double series 
with our normal $M^\la(z)$ is absolutely convergent for
$1/2<\si<3/2$. Indeed,  setting $\rho'=(1/2)+i\ga'$,
$\rho''=(1/2)+i\ga''$ and
$s=\si+i\tau$, the inequalities (\ref{eq:3.6}),
(\ref{eq:3.7}) show that the terms in the double series
are bounded by
$$ 
C(\la,\tau)(|\ga'|+1)^{-\si}(|\ga''|+1)^{-\si}
(|\ga'+\ga''|+1)^{-1+2\si-3/2}.
$$ 
Observing that the number of zeros $\rho=(1/2)\pm
i\ga$ with $n<\ga\le n+1$ is $\cal{O}(\log n)$, the convergence
now follows from a discrete analog of the following simple lemma;
cf.\ Korevaar \cite{Ko07}.
\begin{lemma} \label{lem:3.3}
For real constants $a,\,b,\,c$, the function
$$
\phi(y,v)=(|y|+1)^{-a}(|v|+1)^{-b}(|y+v|+1)^{-c}
$$
is integrable over ${\Bbb R}^2$ if and only if $a+b>1$,  
$a+c>1$, $b+c>1$ and $a+b+c>2$. 
For integrability over ${\Bbb R}_+^2$ the condition $a+b>1$ may
be dropped. 
\end{lemma}
By the lemma, the part of the double sum $\Si^\la_2(s)$ in
(\ref{eq:3.8}) in which $\ga'={\rm Im}\,\rho'$ and $\ga''={\rm
Im}\,\rho''$ have the same sign defines a meromorphic function
for $0<\si<1$ whose only poles occur at the complex zeros of
$\ze(\cdot)$. Thus for a study of its pole-type behavior near the
point $s=1/2$, the sum $\Si^\la_2(s)$ in (\ref{eq:3.8}) may be
reduced to the sum $\Si^\la_3(s)$ in which $\ga'$ and $\ga''$ have
opposite sign. Replacing $\ga''$ by $-\ga''$ and using standard
asymptotics for the Gamma function, it follows that the pole-type
behavior of $\Si^\la_3(s)$ and $\Si^\la(s)$ as $s\searrow 1/2$ is
the same as that of the reduced sum
\begin{equation} \label{eq:3.9}
\Si^\la_4(s) = 2\pi\sum_{\ga'>0,\,\ga''>0}\,
(\ga'\ga'')^{-s+i(\ga'-\ga'')/2}\,M^\la\{1-2s+i(\ga'-\ga'')\}.
\end{equation}
Hence in the study of the PPC under RH, the {\it differences} of
zeta's zeros in, say, the upper half-plane play a key role; cf.\
Montgomery \cite{Mo73}.

Formally, the poles of $\Si^\la(s)$ at the points $s=\rho$ cancel
each other. Under assumption (\ref{eq:2.4}), the function
$\Si^\la_2(s)$ has a meromorphic continuation to the half-plane
$\{\si>1/4\}$, and then there will be real cancellation; see
(\ref{eq:3.1}).

\setcounter{equation}{0} 
\section{Motivation of Conjecture \ref{con:2.2}} \label{sec:4}
As we saw, numerical results make it plausible that the functions
$G_{2r}(s)$ have an analytic continuation to the half-plane
$\{\si>1/4\}$. Hence by Theorem \ref{the:3.1} and
(\ref{eq:3.3}), assuming RH, the function 
\begin{align} \label{eq:4.1}
\Psi^\la(s) & \stackrel{\mathrm{def}}{=} 
2\sum_{0<2r\le\la}\,E(2r/\la)G_{2r}(s) 
+ 4A^E\la \sum_\rho\,\frac{1}{2s-\rho} \notag \\ &
= \Si^\la(s) - D_0(s) - \frac{2R^E(1/2,\la)}{2s-1}+H^\la_0(s) 
\end{align}
will also have an analytic continuation to the half-plane
$\{\si>1/4\}$. This shows that the quadratic pole $1/(2s-1)^2$ of
$D_0(s)$ at $s=1/2$ must be cancelled by a pole of $\Si^\la(s)$ at
$s=1/2$. By (\ref{eq:3.8}) the latter pole is the same as that of
the double sum $\Si^\la_2(s)$, and hence, of
$\Si^\la_4(s)$ in (\ref{eq:3.9}). 

Computation suggests that the quadratic part of the pole at
$s=1/2$ in $\Si^\la_4(s)$ comes from the terms with $\ga''=\ga'$.
Indeed, for $s=(1/2)+\de$ with small $\de>0$, and using the
counting function $N(t)$ of zeta's complex zeros,
\begin{align}
2\pi\sum_{\ga>0}\,\ga^{-1-2\de}M^\la(-2\de) &=
2\pi\la^{-2\de}M^1(-2\de)\int_1^\infty t^{-1-2\de}dN(t)
\notag \\ &=
\frac{1}{4\de^2}-\frac{\log\la+\cal{O}(1)}{2\de}+\cdots.\notag
\end{align}
The residue of the pole of $\Si^\la(s)$ at $s=1/2$ is
$R^E(1/2,\la)$. We know that under assumption (\ref{eq:2.4}),
the residue is $\cal{O}(\la^\eps)$. Independently of (\ref{eq:2.4}),
{\it this may be deduced heuristically} from the fact that $\la$
occurs in the terms of the series for $\Si^\la_2(s)$ only as
$\la^{\rho'+\rho''-2s}$, which is $\cal{O}(\la^\eps)$ for
$s\approx 1/2$. 

We now turn to the likely behavior of $\Psi^\la(s)$ in
(\ref{eq:4.1}) near the line $L=\{\si=1/4\}$.
Since $D_0(s)$ has quadratic poles at the points $s=1/4$ and
$\rho/2$, and no other poles on $L$, cf.\ (\ref{eq:3.2}), we
assume that the (meromorphic continuation of the) double sum 
$\Si^\la_2(s)$ likewise has poles at $1/4$ and the points
$\rho/2$, and nowhere else on $L$. This assumption is plausible because
it is known to be true for $\la\le 2$, when the sum over
$r$ in (\ref{eq:4.1}) is empty, so that the difference
$\Si^\la_2(s)-D_0(s)$ has no poles on $L$ other than first order poles
at the points $\rho/2$. If the {\it heuristic argument} in the preceding
paragraph has general validity, one expects that the coefficients of the
pole terms of $\Si^\la_2(s)$ at $1/4$ and the points $\rho/2$ are
$\cal{O}(\la^{1/2})$, or in any case $\cal{O}(\la^{(1/2)+\eps})$ for
every $\eps>0$. Indeed, the terms in $\Si^\la_2(s)$ contain $\la$ as a
factor $\la^{\rho'+\rho''-2s}=\cal{O}(\la^{(1/2)+\eps})$ for
$\si\searrow 1/4$. 

Hence by (\ref{eq:4.1}), taking coefficient bound $\cal{O}(\la^{1/2})$
for simplicity, the sum $2\sum_{r=1}^N\,E(r/N)G_{2r}(s)$ should behave
like
$$
\frac{\cal{O}(\la^{1/2})}{(4s-1)^2} + \frac{\cal{O}(\la^{1/2})}{4s-1}
$$
near the point $s=1/4$, and like
$$
\frac{\cal{O}(\la^{1/2})}{(2s-\rho)^2} -
\frac{4A^E\la + \cal{O}(\la^{1/2})}{2s-\rho}
$$
near the points $s=\rho/2$. Assuming uniformity here relative to
$\rho$, and taking $\la=2N$, the singular part of the average 
$$
\frac{1}{N}\sum_{r=1}^N\,E(r/N)G_{2r}(s)
$$
for $\si\ge 1/4$ will have the form
\begin{align}  
\frac{1}{N}\sum_{r=1}^N\,G_{2r}(s) &= \frac{\cal{O}(N^{-1/2})}{(4s-1)^2}
+\frac{\cal{O}(N^{-1/2})}{4s-1} +
\cal{O}(N^{-1/2})\sum_{\rho}\,\frac{1}{(2s-\rho)^2} \notag \\ & \quad
-\{4A^E+\cal{O}(N^{-1/2})\}\,\sum_\rho\,\frac{1}{2s-\rho}+H^N(s).\notag
\end{align}
The remainder $H^{N,E}(s)$ will have good boundary behavior as
$\si\searrow 1/4$ and it contains $1/(2N)$ times the remainder 
$H^{\la,E}_0(s)$ from
Theorem \ref{the:3.1} with $\la=2N$. Now taking $E(\nu)$ close to the
function which is equal to $1$ on $[-1,1]$ and $0$ elsewhere, one is led
to Conjecture \ref{con:2.2}.

This conjecture, finally, makes the conjectured Approximation
\ref{approx:2.1} plausible through formal Fourier inversion
(\ref{eq:2.10}).

\setcounter{equation}{0} 
\section{From Approximation \ref{approx:2.1} to Approximation
\ref{approx:1.1} via Approximation \ref{approx:5.2}} \label{sec:5}
After motivating Approximation \ref{approx:2.1} for averages of
functions $\psi_{2r}(x)$, we turn to a corresponding
approximation involving the functions $\theta_{2r}(x)$. For large $N$ and
$x>2N$, cf.\ (\ref{eq:2.5}), 
\begin{align} & \label{eq:5.1}
\frac{1}{N}\sum_{r=1}^N\,\theta_{2r}(x) =
\frac{1}{N}\sum_{r=1}^N\,\psi_{2r}(x) -
\frac{2}{N}\sum_{r=1}^N\,\theta^*_{2r}(x^{1/2}) \notag \\ & \quad
-\frac{1}{N}\sum_{r=1}^N\,4rC_{2r}\log\log x
+\cal{O}(x^{1/3}\log^2 x)+o(N).
\end{align}
According to the Bateman--Horn conjecture
\cite{BH62}, \cite{BH65}, applied to the special case of prime pairs
$(p,\,p^2\pm 2r)$, there should be specific positive constants
$2C^*_{2r}=C^{[2]}_{2r}+C^{[-2]}_{2r}$ such that
\begin{equation} \label{eq:5.2}
\theta^*_{2r}(x)=\{2C^*_{2r}+o(1)\}x\quad\mbox{as}\;\;x\to\infty.
\end{equation}
Here there is no need to study the Bateman--Horn constants in
detail; our only concern will be their mean value (apparently equal
to one).
\begin{conjecture} \label{con:5.1}
For $x\to\infty$ one has
\begin{equation} \label{eq:5.3}
\frac{1}{N}\sum_{r=1}^N\,\theta^*_{2r}(x)
=\{2+o(N^{-1/2})\}x.
\end{equation}
\end{conjecture}

Combining this conjecture with (\ref{eq:5.1}) and Approximation
\ref{approx:2.1}, one obtains the (conjectured)
\begin{approximation} \label{approx:5.2}
For large $N$ and $x$ much larger than $N$, one has
\begin{align} & \label{eq:5.4}
\frac{1}{N}\sum_{r=1}^N\,\{\theta_{2r}(x)-2C_{2r}x\} =
-\{4+\cal{O}(N^{-1/2}\log x)\}\sum_\rho\,\frac{x^\rho}{\rho} 
- \{1+o(1)\}N \notag \\ & 
-\{4+\cal{O}(N^{-1/2}\log x)\}x^{1/2}
-\frac{1}{N}\sum_{r=1}^N\,4rC_{2r}\log\log x + \cal{O}(x^{1/3}\log^2 x).
\end{align}
\end{approximation}
To go from here to Approximation \ref{approx:1.1} we
integrate by parts:
$$
\frac{1}{N}\sum_{r=1}^N\,\pi_{2r}(x)=
\frac{1}{N}\sum_{r=1}^N\,\int_2^x\frac{1}{\log^2 t}\,d\theta_{2r}(t).
$$
In the evaluation it is assumed that contributions due to
derivatives of the $\cal{O}$-terms can be neglected. The sum on
the left of (\ref{eq:5.4}) then becomes the sum on the left of
(\ref{eq:1.5}). Ignoring the $\log\log x$-term for a moment, the
right-hand side of (\ref{eq:5.4}) then gives 
the right-hand side of (\ref{eq:1.5}). The $\log\log x$-term
(with its minus sign) ultimately leads to a contribution
\begin{equation} \label{eq:5.5}
-\frac{1}{N}\sum_{r=1}^N\,4rC_{2r}\int_{2r}^x\frac{1}{\log^2
t}\,d\log\log t\approx -\frac{N}{\log^2 2N}.
\end{equation}
To assess its effect on $\De_N(x)$ in (\ref{eq:1.7}), one still
has to divide by ${\rm li}_2(x^{1/2})\sim 4x^{1/2}/\log^2 x$. The result
$\overline\De_N(x)$ in (\ref{eq:1.8}) will be small when $x$ is much
larger than $N^2$.

In support of Conjecture \ref{con:5.1} we proceed with a
conjecture involving related Dirichlet series
\begin{equation} \label{eq:5.6}
D^*_{2r}(s) = \sum_{p,\,p^2\pm 2r\,{\rm prime}}\,\frac{\log^2
p}{p^{4s}}
=\int_1^\infty x^{-4s}d\theta^*_{2r}(x)\qquad(\si>1/4).
\end{equation}
\begin{conjecture} \label{con:5.3}
For $\si>1/4$ and $N\to \infty$, one has
\begin{equation} \label{eq:5.7}
\frac{1}{N}\sum_{r=1}^N\,D^*_{2r}(s)= \frac{2+o(N^{-1/2})}{4s-1}
+H^N_2(s),
\end{equation}
with an analytic function $H^N_2(s)$ that has good boundary
behavior as $\si\searrow 1/4$.  
\end{conjecture}

The arguments supporting Conjecture \ref{con:5.3} are similar to
those given for Conjecture \ref{con:2.2}. For $1/4<\si<1/2$ one
may write
\begin{align} & \label{eq:5.8}
\Phi^\la_{1,2}(s) \stackrel{\mathrm{def}}{=} D^*_0(s)+
2\sum_{0<2r\le\la}\,E(2r/\la)D^*_{2r}(s) \notag \\ &=
\frac{2A^E\la}{4s-1}+\frac{\ze'(2s)}{\ze(2s)}\,J(s,s)+
\sum_\rho\,\frac{1}{2}\Ga\{(\rho/2)-s\}J(\rho/2,s)+H^\la_3(s),
\end{align}
where $D^*_0(s)=\sum_p\,(\log^2 p)/p^{4s}$ and $H^\la_3(s)$ is
analytic for $1/4\le\si<1/2$. The functions $J(s,s)$ and
$J(\rho/2,s)$ are analytic for $1/4\le\si<1/2$; cf.\
(\ref{eq:8.4}) in the Appendix. Hence, formally the poles at the
points $s=\rho/2$ in the combination of $J$-terms in
(\ref{eq:5.8}) will cancel each other. However, one constituent
of $J(\rho/2,s)$ is an infinite series. It leads to a repeated
series $\Si^\la_{2,2}(s)$ when it is substituted into the sum
over $\rho$ in (\ref{eq:5.8}):
\begin{align} 
\Si^\la_{2,2}(s) &= 
\frac{1}{2}\sum_{\rho}\,\Ga\{(\rho/2)-s\}\sum_{\rho'}\,\Ga(\rho'-s)
\notag \\ & \qquad \times
M^\la\{\rho'+(\rho/2)-2s\}\cos\{\pi(\rho'-\rho/2)/2\}.\notag
\end{align}
This series is absolutely convergent only for $3/8<\si<1/2$; for
$1/4<\si\le 3/8$ the sum over $\rho=(1/2)+i\ga$ has to be
interpreted as a limit of partial sums $\sum_{|\ga|\le B}$ as
$B\to\infty$. 

In view of the similarity of the Hardy--Littlewood conjecture and
our case of the Bateman--Horn conjecture, it is reasonable to
suppose that the differences
$$
G^*_{2r}(s)=D^*_{2r}(s)-\frac{2C^*_{2r}}{4s-1}
$$
have an analytic continuation to the half-plane $\{\si\ge
1/4\}$. If that is correct, the combination of the $J$-terms in
(\ref{eq:5.8}) truly has no poles at the points $s=\rho/2$.
The repeated sum $\Si^\la_{2,2}(s)$ then would have an analytic
continuation to the strip $1/4\le\si<1/2$, except for a pole at
$s=1/4$. The quadratic pole
$1/(4s-1)^2$ of $D^*_0(s)$ at $s=1/4$ would be cancelled by the
quadratic part of the pole of $\Si^\la_{2,2}(s)$ there. Finally,
the residue of the pole of $\Si^\la_{2,2}(s)$ at $s=1/4$ would be
$\cal{O}(\la^{(1/4)+\eps})$ by heuristics as in Section \ref{sec:4}.
Hence by (\ref{eq:5.8}), the residue at $s=1/4$ of 
$$
\frac{1}{N}\sum_{r=1}^N\,E(r/N)D^*_{2r}(s)\quad\mbox{would
be}\quad \frac{1}{2}A^E+o(N^{-1/2}),
$$
thus leading to (\ref{eq:5.7}) when $E(\nu)$ is taken close to $1$ on
$(-1,1)$.

The proof of (\ref{eq:5.8}) is similar to that of (\ref{eq:3.1})
described in the Appendix. Here, one would start with the
integral obtained from (\ref{eq:8.2}) through replacement of one
of the quotients $\ze'(\cdot)/\ze(\cdot)$ by
$\ze'(2\,\cdot)/\ze(2\,\cdot)$.

\newpage
\begin{center}
{\bf PART II. NUMERICAL RESULTS AND GRAPHS}
\end{center}

\setcounter{equation}{0} 
\section{Comparing averages of functions $\om_{2r}(x)$ with
$\om(x)$}
\label{sec:6}
Which of the two terms on the right-hand side
of (\ref{eq:1.5}), in the conjectured Approximation \ref{approx:1.1}, is
larger? One may write
\begin{align}
\sum_\rho\,\rho\,{\rm li}_2(x^\rho) &= \frac{x^{1/2}}{\log^2 x}
\sum_\rho\,\frac{x^{\rho-1/2}}{\rho}
+\cal{O}\bigg(\frac{x^{1/2}}{\log^3 x}\bigg),
\notag \\ 
{\rm li}_2(x^{1/2}) &= 4\frac{x^{1/2}}{\log^2 x}
+\cal{O}\bigg(\frac{x^{1/2}}{\log^3 x}\bigg).\notag
\end{align}
Hence if we set
\begin{equation} \label{eq:6.1}
T(x) \stackrel{\mathrm{def}}{=}
\sum_\rho\,\frac{x^{\rho-1/2}}{\rho}=\sum_{\ga>0}\,
\frac{\cos(\ga\log x)+2\ga\sin(\ga\log x)}{\ga^2+1/4},
\end{equation}
relation (\ref{eq:1.5}) takes the form 
\begin{align} \label{eq:6.2}
Q_N(x) &\stackrel{\mathrm{def}}{=}
\frac{\sum_{r=1}^N\,\om_{2r}(x)}{N{\rm
li}_2(x^{1/2})} \notag \\ &= -[\{1+\cal{O}(N^{-1/2}\log
x)\}T(x)+1+\cal{O}(N^{-1/2}\log x)].
\end{align}
  
It is interesting that Riemann's formula (\ref{eq:1.9}) 
leads to a combination similar to the right-hand side of
(\ref{eq:6.2}). Indeed, assuming RH one may write
\begin{equation} \label{eq:6.3}
\frac{2\om(x)}{{\rm
li}(x^{1/2})}=-\big[\{1+\cal{O}(1/\log
x)\}T(x)+1+\cal{O}(1/\log x)\big].
\end{equation}
Littlewood's work \cite{Li14}, cf.\ Ingham \cite{In90}, implies that the
function
$T(x)$ oscillates unboundedly. More precisely, he showed that there are
constants $c,\,c'>0$ and arbitrarily large $x,\,x'$ such that
$$
T(x)<-c\log\log\log x, \quad T(x')>c'\log\log\log x'.
$$
However, $\pi(x)$ becomes larger than ${\rm li}(x)$, that is,
$\om(x)>0$, only for certain very large $x$. The first such
number is associated with the name of Skewes; cf.\ te
Riele \cite{tR87}, and Bays and Hudson \cite{BH00}. Under RH one has
$T(x)=\cal{O}(\log^2 x)$ and Kotnik \cite{Kt08} made it
plausible that $T(x)=\cal{O}(\log x)$. He also graphed the function
$\om(x)(\log x)/x^{1/2}$, cf.\ (\ref{eq:6.3}), for
$x\le 10^{14}$. On a logarithmic scale, his Figure 1 shows rapid
oscillations of amplitude greater than $1/2$.

The Skewes story seems to have no analog for prime pairs, cf.\
Brent \cite{Br75}. Here we focus on the case of twin primes. 
Nicely \cite{Ni08} has counted prime twins up
to $x=10^{16}$. His table uses steps $10^k$ from
$1\cdot 10^k$ through $9\cdot 10^k$ for
$k=1,2,\cdots,12$. From there on the steps are $10^{12}$.
Nicely's table shows that for $x$ going to $10^{16}$, the
quantity $|\om_2(x)|$ often becomes a good deal larger than 
${\rm li}_2(x^{1/2})$.
His table implies $16$ {\it sign changes} of
$\om_2(x)$ [which is minus his entry $\de_2(x)$]. The first occurs
between $10^6$ and $2\cdot 10^6$, the last between $7.5\cdot
10^{13}$ and $7.6\cdot 10^{13}$. Although
$\om_2(x)$ oscillates, it then remains positive until the end of
Nicely's table. 

In a preprint on a `Skewes number for twin primes', Marek Wolf
\cite{Woxx} analyzed the sign changes in $\om_2(x)$ up to
$2^{42}\approx 4.4\cdot 10^{12}$. He found the first one at the
twin with $p=1369391$. A table in his preprint lists the number
of sign changes up to $2^k$ for $k=22,23,\cdots,42$. Wolf found
$90355$ sign changes up to $2^{42}$. He found none between
$2^{22}$ and $2^{25}$, none between $2^{28}$ and $2^{31}$, and
none between $2^{37}$ and $2^{39}$.

In our range of $x$, the values of $|T(x)|$ are smaller than
one. In particular 
\begin{align} & \label{eq:6.4}
T(10^6)\approx 0.41156,\quad T(10^8)\approx 0.17554,\notag \\ &
T(10^{10})\approx -0.42122,\quad T(10^{12})\approx -0.04014.
\end{align}
These values were computed with the aid of von Mangoldt's formula 
(\ref{eq:1.10}), by which (for $x>1$ and $x$ not a prime power)
\begin{equation} \label{eq:6.5}
T(x)=x^{-1/2}\big\{x-\psi(x)-\log(2\pi)-(1/2)\log\,(1-x^{-2})\big\}.
\end{equation}

The function $\psi(x)=\sum_{p^m\le x}\log p$ was computed by
summing the values of $\lfloor \log_px\rfloor\log p$
for all the primes $p\le x$ (generated with the sieve of Eratosthenes).
Here, $\lfloor \log_px\rfloor$ is the exponent of $p$ in the highest power of $p$ not exceeding $x$.
The values of $T(x)$ given in (\ref{eq:6.4})
were computed with an accuracy of at least 5 decimal digits.
We were using Fortran double precision floating point arithmetic
which works with an accuracy of about 15 decimal digits, but precision
is lost as $x$ grows when (\ref{eq:6.5}) is used to compute $T(x)$.
To illustrate this, we found that $\psi(10^{12})=1000000040136.76$, so that
in the difference $10^{12}-\psi(10^{12})=-40136.76$ only about seven digits
are still correct and $T(10^{12})=-0.04013860$.

Alternative computations based on formula (\ref{eq:6.1}) and the first 
two million values of $\ga$ gave the values $T(10^6)\approx0.41276$,
$T(10^8)\approx0.17469$, $T(10^{10})\approx-0.41944$, and
$T(10^{12})\approx-0.04010$, i.e., an accuracy of only about 3 decimal digits.

\setcounter{equation}{0} 
\section{Testing the conjectured Approximation \ref{approx:1.1}}
\label{sec:7}
In the following we will consider the aggregate
\begin{equation} \label{eq:7.1}
\Pi_N(x)\stackrel{\mathrm{def}}{=}\pi_2(x)+\pi_4(x)+\cdots+
\pi_{2N}(x) 
\end{equation}
for certain large values of $N$ and $x$. Setting
\begin{equation} \label{eq:7.2}
S_N=C_2+C_4+\cdots+C_{2N}, 
\end{equation}
cf.\ (\ref{eq:3.4}), we compare $\Pi_N(x)$ with $2S_N\,{\rm
li}_2(x)=(S_N/C_2)L_2(x)$. In view of the conjectured Approximation
\ref{approx:1.1}, the difference is divided by $N\,{\rm li}_2(x^{1/2})$
to obtain the quotient 
\begin{equation} \label{eq:7.3}
\frac{\Pi_N(x)-(S_N/C_2)L_2(x)}{N\,{\rm
li}_2(x^{1/2})}=\frac{\sum_{r=1}^N\,\om_{2r}(x)}{N\,{\rm
li}_2(x^{1/2})}=Q_N(x);
\end{equation}
cf.\ (\ref{eq:6.2}). For large $N$ the quotient
should have the form
$$
-\{1+\cal{O}(N^{-1/2}\log x)\}T(x)-\{1+\cal{O}(N^{-1/2}\log x)\}
$$
Ignoring the $\cal{O}$-terms, we 
will compare $Q_N(x)$ with $-T(x)-1$, setting
\begin{equation} \label{eq:7.4}
Q_N(x)+T(x)+1=\De_N(x).
\end{equation}

Tables 2, 3 give results for $x=10^6,\,10^8,\,10^{10},\,10^{12}$. The
values $S_N/C_2$ were obtained by computing $C_{2r}/C_2$ from
(\ref{eq:1.3}) and adding. For the values of 
\begin{table} \label{table:2}
\begin{tabular}{rrrrrrr} 
$2N$ & $S_N/C_2$ & $\Pi_N(10^6)$ & $\De_N(10^6)$ & $\Pi_N(10^8)$ &
$\De_N(10^8)$ & \\
       &             &         &          &          &            & \\
$100$  &  73.6377551 &  605087 & +0.09722 &  32417440 & -0.08872  & \\
$200$  & 149.3252708 & 1226667 & -0.02199 &  65739481 & +0.03162  & \\
$300$  & 225.4407734 & 1851433 & -0.12785 &  99245855 & -0.09833  & \\
$400$  & 300.3132204 & 2465581 & -0.23344 & 132202659 & -0.23013  & \\ 
$500$  & 376.0636735 & 3086695 & -0.32860 & 165551273 & -0.18188  & \\
$600$  & 452.4693143 & 3714028 & -0.31371 & 199186203 & -0.19507  & \\
$700$  & 527.3827110 & 4328507 & -0.34805 & 232164862 & -0.18926  & \\
$800$  & 603.4536365 & 4951873 & -0.42140 & 265651152 & -0.21737  & \\
$900$  & 679.4011178 & 5574196 & -0.48004 & 299079601 & -0.28690  & \\
$1000$ & 754.4223630 & 6188960 & -0.52230 & 332105577 & -0.27582  & \\ 
$2000$ & 1511.5853400 & 12391586 & -0.78001 & 665435604 & -0.16751 & \\
$3000$ & 2269.6853566 & 18597363 & -0.95390 & 999175096 & -0.14446 & \\
$4000$ & 3026.0445409 & 24783891 & -1.11135 & 1332114654 & -0.23565 & \\
$5000$ & 3783.8474197 & 30975067 & -1.28953 & 1665693721 & -0.28111 & \\
       &             &         &          &           &           &
\end{tabular}
\caption{Values of $S_N/C_2$, $\Pi_N(10^6)$, $\De_N(10^6)$ and
$\Pi_N(10^8)$, $\De_N(10^8)$}
\end{table}
$\Pi_N(x)$ we added columns of numbers $\pi_{2r}(x)$. We
next computed $Q_N(x)$ from (\ref{eq:7.3}). Here we used the
approximations
\begin{align} &
L_2(10^6)\approx 8 248.0297,\quad
L_2(10^8)\approx 440 367.7942, \notag \\ &
L_2(10^{10})\approx 27 411 416.53,\quad 
L_2(10^{12})\approx 1 870 559 866.82 \notag
\end{align} 
and 
\begin{align} &
{\rm li}_2(10^3)\approx 34.6851,\quad {\rm
li}_2(10^4)\approx 162.2412, \notag \\ &
{\rm li}_2(10^5)\approx 945.75959, \quad
{\rm li}_2(10^6)\approx 6 246.9757 \notag.
\end{align}
The table entries
$\De_N(x)$ are based on (\ref{eq:7.4}) and the approximations for
$T(x)$ in (6.4). 
\begin{table} \label{table:3}
\begin{tabular}{rrrrrr} 
$2N$ & $\Pi_N(10^{10})$ & $\De_N(10^{10})$ & $\Pi_N(10^{12})$ &
$\De_N(10^{12})$ & \\
       &              &           &                &           & \\
$100$  &  2018498733  & +0.23101  &  137743459486  & -0.22449  & \\
$200$  &  4093181354  & +0.19981  &  279320931774  & -0.52374  & \\
$300$  &  6179575427  & +0.04646  &  421698995095  & -0.60678  & \\
$400$  &  8231900717  & -0.00307  &  561752066806  & -0.47345  & \\ 
$500$  & 10308323520  & +0.09461  &  703447298670  & -0.52336 & \\
$600$  & 12402663153  & +0.00891  &  846368266787  & -0.46665  & \\
$700$  & 14456137134  & +0.06512  &  986498011024  & -0.37686  & \\
$800$  & 16541312091  & +0.03187  & 1128792535379  & -0.48827  & \\
$900$  & 18623097684  & -0.00710  & 1270856645797  & -0.39850  & \\
$1000$ & 20679532323  & +0.04311  & 1411187901897  & -0.41454  & \\
$2000$ & 41434008965  & -0.14700  & 2827502930522  & -0.31142  & \\
$3000$ & 62214267139  & -0.14273  & 4245571295213  & -0.21865  & \\
$4000$ & 82946817735  & -0.13473  & 5660383932743  & -0.12392  & \\
$5000$ & 103718886923 & -0.15324  & 7077896171945  & -0.12569  & \\ 
       &              &           &                &           &
\end{tabular}
\caption{Values of $\Pi_N(10^{10})$, $\De_N(10^{10})$ and
$\Pi_N(10^{12})$,
$\De_N(10^{12})$}
\end{table}

In Figures 1--4 we show plots of $\De_N(x)$ as a function of $N$
($50\le 2N\le 5000$), for $x=10^6$, $10^8$, $10^{10}$, and $10^{12}$.
We have omitted the function values for $2\le 2N\le48$ since they very much dominate
(and are atypical for) the other function values.
In Figures 1 and 2 we compare $\De_N(x)$ with the function
$\overline\De_N(x)=-\frac{2N\log^2 x}{8x^{1/2}\log^2 2N}$
as defined in (\ref{eq:1.8}).

We have made, but not given here, plots of
$\De_N(x)$ for several other values of $x$.
E.g., for $x=10^{11}$ and $2N=1000, 2000, 3000, 4000, 5000$, we found:
$\De_N(x)=-0.229, -0.072, -0.034, +0.004$, and $-0.034$, respectively
(compare these values with the corresponding values for $x=10^{10}$
and $x=10^{12}$ in Table 3).

Figures 5 and 6 show plots of $\De_N(x)$ as a function of $x$
($6\le\log_{10} x\le12$), for $N=400$ and $N=2500$, respectively.
The plots have been constructed by connecting the
values of $\De_N(x)$ for $x=10^6$ and for $x=i\times10^j$, $j=6,7,\dots,11$
and $i=1,2,\dots,10$ by straight lines.
The different behaviour of the plots of $\De_{400}(x)$ and $\De_{2500}(x)$
may reflect the influence of the unknown $\cal{O}(N^{-1/2}\log x)$--terms,
which were neglected in the derivation of the error function $\De_N(x)$
from (\ref{eq:1.5}).

\newpage
\begin{center}
\includegraphics[angle=-90,totalheight=8.5cm]{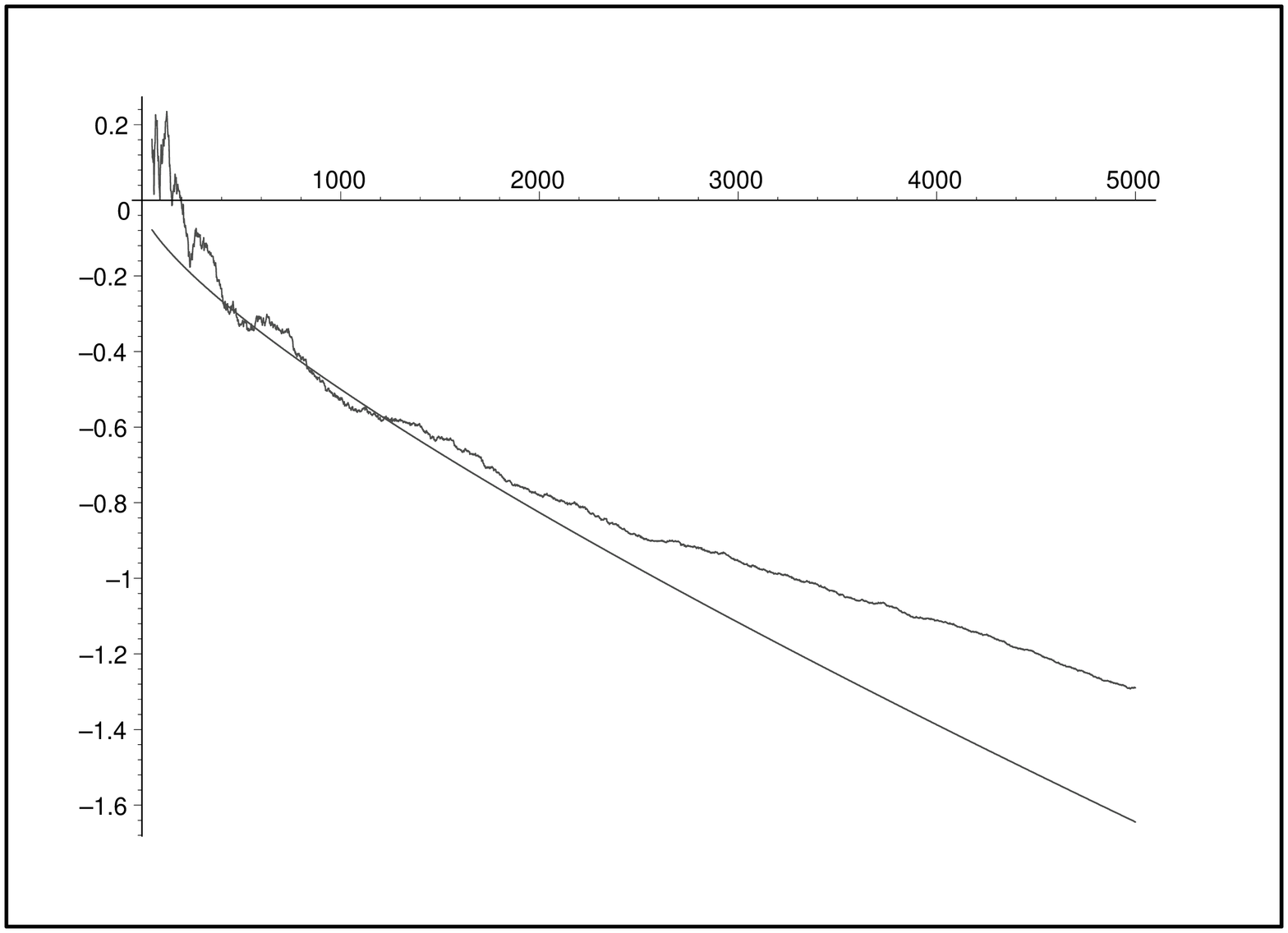}
\\{\sc Figure 1}. $\De_N(10^6)$ compared with $\overline\De_N(10^6)$ for $50\le2N\le5000$
\end{center}
\begin{center}
\includegraphics[angle=-90,totalheight=8.5cm]{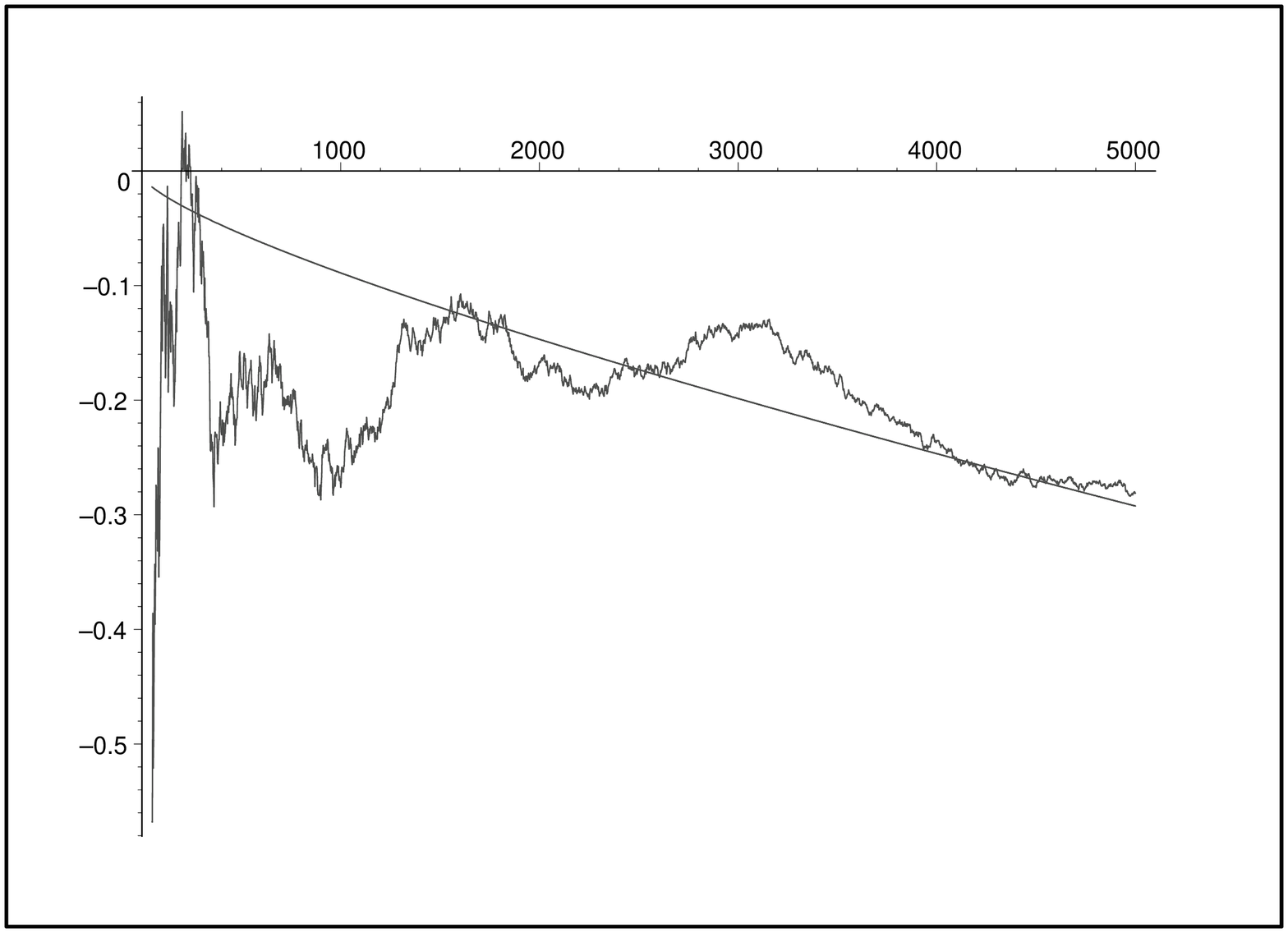}
\\{\sc Figure 2}. $\De_N(10^8)$ compared with $\overline\De_N(10^8)$ for $50\le2N\le5000$
\end{center}
\begin{center}
\includegraphics[angle=-90,totalheight=8.5cm]{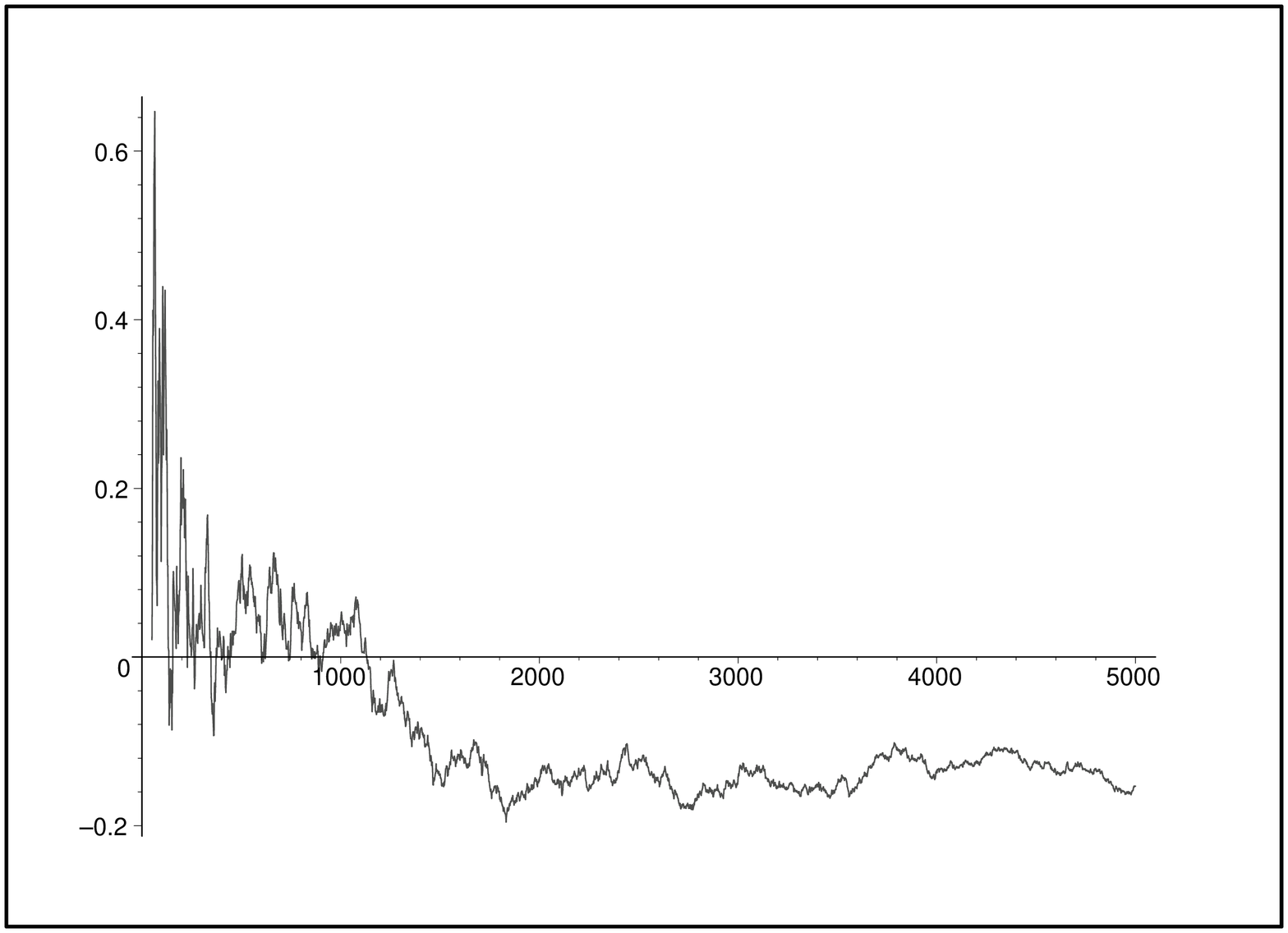}
\\{\sc Figure 3}. $\De_N(10^{10})$ for $50\le2N\le5000$
\end{center}
\begin{center}
\includegraphics[angle=-90,totalheight=8.5cm]{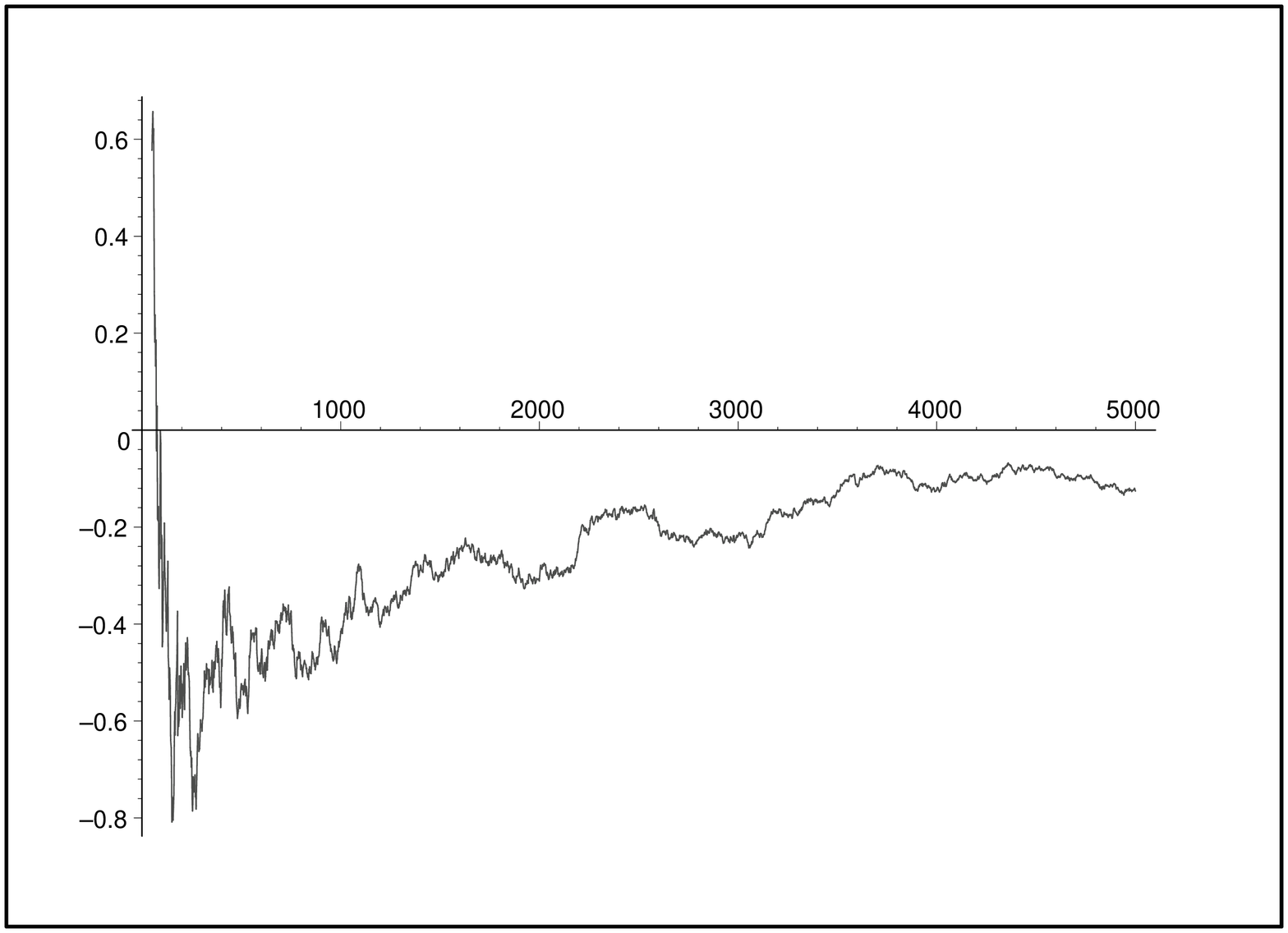}
\\{\sc Figure 4}. $\De_N(10^{12})$ for $50\le2N\le5000$
\end{center}
\begin{center}
\includegraphics[angle=-90,totalheight=8.5cm]{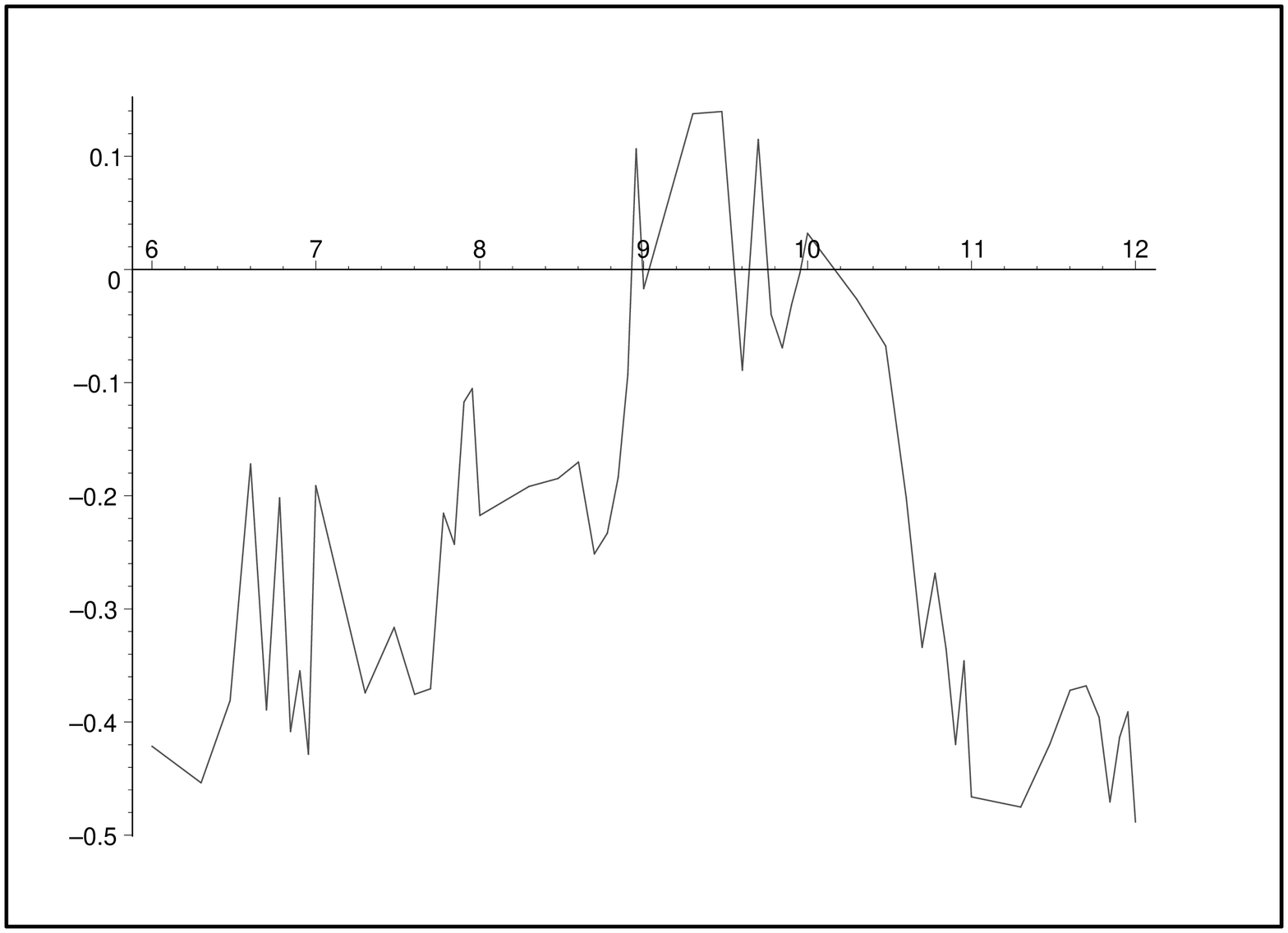}
\\{\sc Figure 5}. $\De_{400}(x)$ for $6\le \log_{10} x\le12$
\end{center}
\begin{center}
\includegraphics[angle=-90,totalheight=8.5cm]{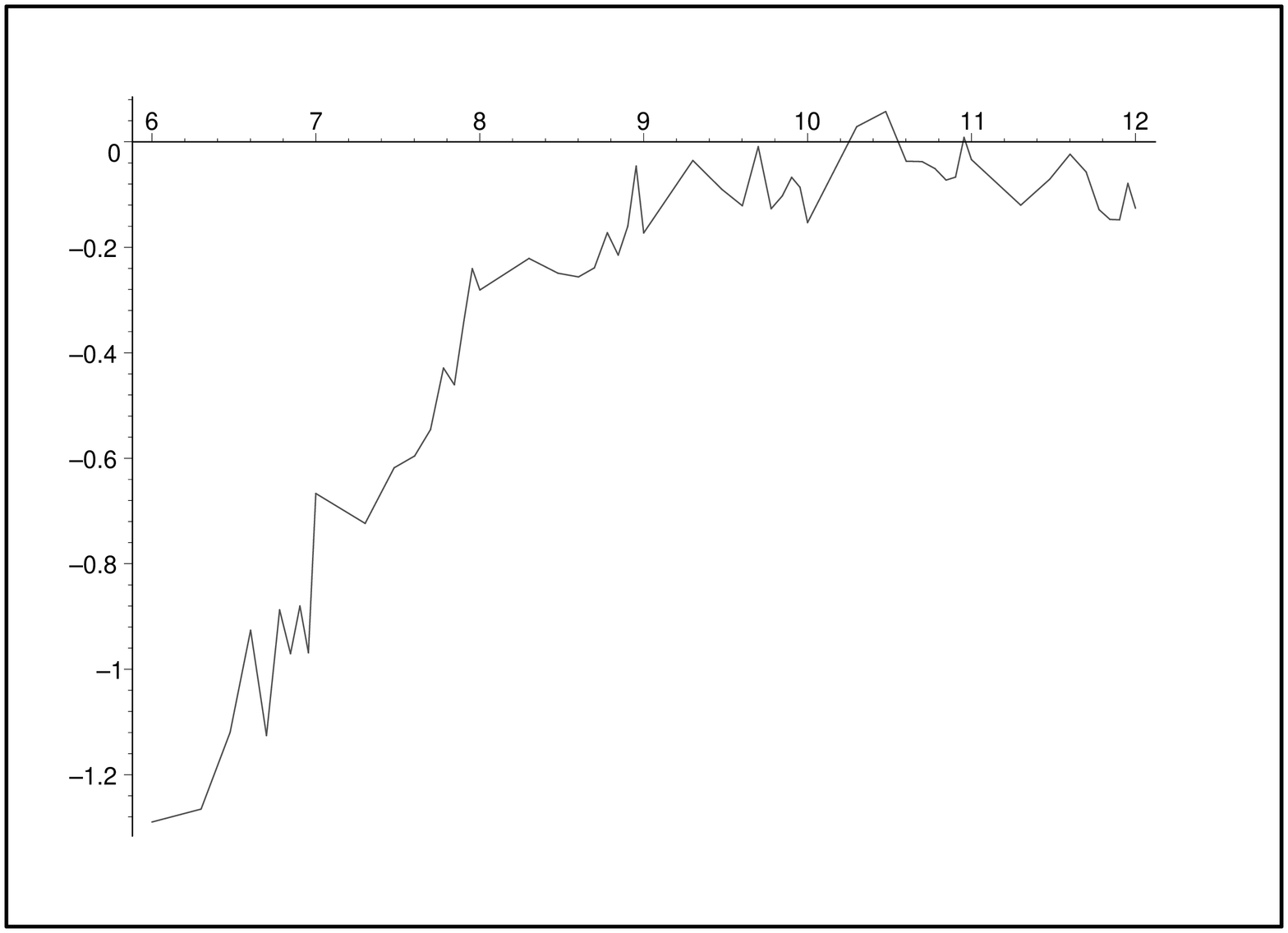}
\\{\sc Figure 6}. $\De_{2500}(x)$ for $6\le \log_{10} x\le12$
\end{center}

\newpage
\begin{center}
{\bf APPENDIX}
\end{center}

\setcounter{equation}{0} 
\section{Outline of the proof of Theorem \ref{the:3.1}}
\label{sec:8} 
Extending an idea that goes back to Arenstorf
\cite{Ar04}, cf.\ \cite{Ko07}, one is led to a representation for
$E\{(\al-\be)/\la\}$ by an absolutely convergent repeated complex
integral in which $\al>0$ and $\be>0$ occur separately:
\begin{align} \label{eq:8.1}
E\{(\al-\be)/\la\} & = \frac{1}{(2\pi
i)^2}\int_{L(c,B)}\Ga(z)\al^{-z}dz
\int_{L(c,B)}\Ga(w)\be^{-w}\notag \\ & \qquad\qquad\quad\times
M^\la(z+w)\cos\{\pi(z-w)/2\}\,dw.
\end{align}
Here the path $L(c,B)=L(c_1,c_2,B)$ in the $z=x+iy$ plane is taken of
the form
$$ 
L(c,B)=\left\{\begin{array}{lllll}
\mbox{$\quad$the half-line} & \mbox{$\{x=c_1,\,-\infty<y\le
-B\}$}\\
\mbox{$+\;$the segment} & \mbox{$\{c_1\le x\le
c_2,\,y=-B\}$}\\
\mbox{$+\;$the segment} & \mbox{$\{x=c_2,\,-B\le y\le B\}$}\\
\mbox{$+\;$the segment} & \mbox{$\{c_2\ge x\ge
c_1,\,y=B\}$}\\
\mbox{$+\;$the half-line} & \mbox{$\{x=c_1,\,B\le
y<\infty\}$,}  
\end{array}\right.
$$
and similarly for the $w=u+iv$ plane. For $-1/2<c_1<0<c_2<1/2$, say, 
and arbitrary $B>0$, the absolute convergence of the
repeated integral in (\ref{eq:8.1}) follows from (\ref{eq:3.6}),
(\ref{eq:3.7}) and Lemma \ref{lem:3.3}.

For the verification of formula (\ref{eq:8.1}) one may write $\cos\al t$
as a complex (inverse) Mellin integral involving $\Ga(z)$, and $\cos\be
t$ as such an integral involving $\Ga(w)$. Multiplying the two, doing the
same with sines and adding, one obtains a repeated complex integral for
$\cos (\al-\be)t$. It is multiplied by $\hat E^\la(t)$;
integration over $0<t<\infty$ and use of (\ref{eq:3.5}) then gives
the result.     

Formula (\ref{eq:8.1}) leads to the following integral for $\Phi^\la(s)$
in (\ref{eq:3.1}), modulo a function $H^\la(s)$ that turns out to be
analytic for $\si>0$:
\begin{align}  \label{eq:8.2}
\quad \Phi^\la(s) &= \frac{1}{(2\pi i)^2}
\int_{L(c,B)}\Ga(z)\,\frac{\ze'(z+s)}{\ze(z+s)}\,dz
\int_{L(c,B)}\Ga(w)\,\frac{\ze'(w+s)}{\ze(w+s)}\notag \\ &
\qquad\qquad\quad\times M^\la(z+w)\cos\{\pi(z-w)/2\}dw + H^\la(s).
\end{align}
For verification one introduces the Dirichlet series $\sum
\La(k)k^{\cdots}$ and $\sum \La(l)l^{\cdots}$ for the quotients
$\ze'/\ze$. One then integrates term by term, initially taking
$\si>1+|c_1|$. The result
$$
\sum_{k,l=1}^\infty\,\frac{\La(k)\La(l)}{k^s l^s}\,E\{(k-l)/\la\}
$$
differs from $\Phi^\la(s)$ in (\ref{eq:3.1}) by
\begin{align}
H^\la(s) &= 2 \sum_{0<2r\le\la}\,\sum_{n=1}^\infty\,\La(n)\La(n+2r)
\bigg\{\frac{1}{n^{2s}}-\frac{1}{n^s(n+2r)^s}\bigg\}E(2r/\la)
\notag \\ & \quad
-2 \sum_{0<2r-1\le\la}\,\sum_{n=1}^\infty\,\frac{\La(n)\La(n+2r-1)}
{n^s(n+2r-1)^s}\,E\{(2r-1)/\la\}.\notag
\end{align} 
Thus $H^\la(s)$ is analytic for $\si>0$. One may verify that it
is the Mellin transform of a function
$h^\la(x)$ which is $\cal{O}\{\la(1+\la^2/ x)\log^2 x\}$ for $x>\la$.

Analytic continuation shows that under RH, one may take paths
$L(c,B)$ in (\ref{eq:8.2}) with $c_1=-\eta$ and $c_2=(1/2)-\eta$, where 
$0<\eta<1/2$. Thus the integral representation may be used for
$s=\si+i\tau$ with $\si>(1/2)+\eta$ and $|\tau|<B$; cf.\ \cite{Ko07}. 
Additionally requiring $\si<1$, we now move the paths $L(c,B)$
across the poles at the points $1-s$, $0$ and $\rho-s$ to lines
$L(d)$, given by $x$ or $u$ equal to $d=-(1/2)+\eta$. The moves may
be justified by Cauchy's theorem and the estimates in Lemma
\ref{lem:3.3}. On the relevant vertical lines, $(\ze'/\ze)(Z)$ only grows
logarithmically in $Y$, and auxiliary horizontal segments can be
suitably chosen between zeta's complex zeros.

First moving the $w$-path one obtains a new repeated integral,
along with a single `residue-integral'. It is convenient to write the
latter in the form  
\begin{equation} \label{eq:8.3}
 \frac{1}{2\pi i}\int_{L(c,B)}
\Ga(z)\,\frac{\ze'(z+s)}{\ze(z+s)}\,J(z+s,s)dz,
\end{equation}
where by the residue theorem
\begin{align} \label{eq:8.4}
J(z+s,s) &=
-\Ga(1-s)M^\la(z+1-s)\cos\{\pi(z+s-1)/2\}\notag \\ &
\quad +
\frac{\ze'(s)}{\ze(s)}\,M^\la(z)\cos(\pi z/2) \\ &
\quad +\sum_\rho\,\Ga(\rho-s)M^\la(z+\rho-s)
\cos\{\pi(z+s-\rho)/2\}.\notag
\end{align}

Next move the $z$-path $L(c,B)$ in the new repeated integral
and the $z$-path in the single integral to the line $L(d)$. Thus
we obtain another repeated integral, now involving two paths
$L(d)$, and a single integral with path $L(d)$, where
$d=-(1/2)+\eta$. Varying $\eta\in(0,1/2)$, one sees that the new
integrals represent analytic functions for $-1/2<\si<1$. The
operation on the repeated integral produces a harmless residue,
namely, another copy of the single integral with path
$L(d)$. However, the operation on the single integral yields the
following residue:
\begin{equation} \label{eq:8.5}
-\Ga(1-s)J(1,s)+\frac{\ze'(s)}{\ze(s)}\,J(s,s) +
\sum_{\rho'}\,\Ga(\rho'-s)J(\rho',s).
\end{equation}
Working out this residue with the aid of (\ref{eq:8.4}) one
obtains nine terms. Four of these supply the sum $\Si^\la(s)$ of
(\ref{eq:3.8}) in (\ref{eq:3.1}). The remaining five terms
combine into the sum
\begin{align} \label{eq:8.6}
V^\la(s) &\stackrel{\mathrm{def}}{=}
\Ga^2(1-s)M^\la(2-2s)-2\Ga(1-s)\,\frac{\ze'(s)}{\ze(s)}\,
M^\la(1-s)\sin(\pi s/2) \notag \\ & \quad
-2\Ga(1-s)\sum_\rho\,\Ga(\rho-s)M^\la(1+\rho-2s)\sin(\pi\rho/2).
\end{align}
Here, the apparent poles at the points $s=\rho$ cancel each
other. The first term provides the important
pole-term at the point $s=1/2$ in (\ref{eq:3.1}). Indeed, by the
pole-type behavior of $M^\la(Z)$ at the point $Z=1$ (Section
\ref{sec:2}),
\begin{equation} \label{eq:8.7}
\Ga^2(1-s)M^\la(2-2s)=\frac{2A^E\la}{2s-1}+H^\la_4(s),
\end{equation}
where $H^\la_4(s)$ is analytic for $-1/2<\si<1$. The
final term in (\ref{eq:8.6}) generates simple poles at the
points $s=\rho/2$. A short computation shows that the
residues at those poles are all equal to $-2A^E\la$, thus
leading to the term $-4A^E\la\sum_\rho\,1/(2s-\rho)$
in (\ref{eq:3.1}).

To round out the proof of Theorem \ref{the:3.1} we evaluate the inverse
Mellin transform of $V^\la(s)$ in (\ref{eq:8.6}). Taking $c=(1/2)+\eta$
and $x>\la$, and moving the path to the left, one finds that
$$ 
\frac{1}{2\pi i}\int_{L(c)} V^\la(s)\,x^{2s}\,\frac{ds}{s} 
= 2A^E\la\,x - 4A^E\la\sum_\rho\frac{x^\rho}{\rho}+g(x,\la).
$$
Here for large $\la$ and $x/\la\to\infty$, 
$$
g(x,\la)=-\{1+o(1)\}\la^2\int_0^1 E(\nu)\nu\,d\nu\approx -\la^2/2
$$
when $E(\nu)$ is close to $1$ on $(-1,1)$. 

By our hypothesis the double sum $\Si^\la_2(s)$ in (\ref{eq:3.8})
generates a pole at $s=1/2$ with residue $R^E(1/2,\la)$ given by
(\ref{eq:3.3}). Through Mellin inversion this becomes $2R^E(1/2,\la)x$;
when added to $2A^E\la x$, it gives the principal part 
$$
2\sum_{0<2r\le\la}\,E(2r/\la)\cdot 2C_{2r}x\quad\mbox{of}\quad
2\sum_{0<2r\le\la}\,E(2r/\la)\psi_{2r}(x).
$$
Other plausible contributions of $\Si^\la_2(s)$ have been discussed in
Section \ref{sec:4}.

\noindent{{\it KdV Institute of Mathematics, University of 
Amsterdam,} \\
Plantage Muidergracht 24, 1018 TV Amsterdam, The Netherlands}

\noindent{{\it CWI: Centrum Wiskunde en Informatica,} \\
Science Park 123, 1098 XG Amsterdam, The Netherlands}

\noindent{\it E-mail addresses}: {\tt J.Korevaar@uva.nl}, {\tt
Herman.te.Riele@cwi.nl}

\enddocument